\documentclass[11pt,reqno]{amsart}
\usepackage[a4paper,margin=28mm]{geometry}
\usepackage{amsmath,amssymb,mathtools}
\usepackage{xcolor}
\usepackage{hyperref}
\hypersetup{colorlinks=true,linkcolor=blue!40!black,
  citecolor=blue!40!black,urlcolor=blue!40!black,
  pdfauthor={Liding Huang, Jiaogen Zhang, Xi Zhang},
  pdftitle={A second-order estimate linear in the gradient for complex Hessian equations}}
\allowdisplaybreaks[2]
\numberwithin{equation}{section}
\newtheorem{theorem}{Theorem}[section]
\newtheorem{proposition}[theorem]{Proposition}
\newtheorem{lemma}[theorem]{Lemma}
\newtheorem{corollary}[theorem]{Corollary}
\theoremstyle{remark}
\newtheorem{remark}[theorem]{Remark}
\DeclareMathOperator{\osc}{osc}
\DeclareMathOperator{\tr}{tr}
\newcommand{\ddbar}{\sqrt{-1}\partial\bar\partial}
\newcommand{\Lin}{\mathcal L}
\newcommand{\norm}[1]{\left\lVert#1\right\rVert}
\newcommand{\data}{\mathcal D}

\title[A second-order estimate for complex Hessian equations]
{A second-order estimate linear in the gradient\newline
for complex Hessian equations}

\author[L. Huang]{Liding Huang}
\address{School of Mathematical Sciences, Xiamen University,
  Xiamen 361005, P. R. China}
\email{huangliding@xmu.edu.cn}

\author[J. Zhang]{Jiaogen Zhang}
\address{ School of Mathematics, Hefei University of Technology,
  Hefei 230009, P. R. China}
\email{zjgmath@hfut.edu.cn}

\author[X. Zhang]{Xi Zhang}
\address{School of Mathematics and Statistics,
  Nanjing University of Science and Technology, Nanjing 210094, China}
\email{mathzx@njust.edu.cn}

\date{}
\subjclass[2020]{Primary 32W20; Secondary 35J60, 35B45.}
\keywords{Complex Hessian equation, second-order estimate,
maximum principle, stability}

\begin{document}
\begin{abstract}
We use the maximum principle to prove a second-order estimate with linear dependence on the gradient for the nondegenerate complex Hessian equation on a compact K\"ahler manifold. The proof combines $C^0$ compactness and a refined concavity inequality with the use of nearby smooth admissible functions. A standard interpolation argument then yields the known uniform gradient and complex Hessian bounds. 
\end{abstract}
\maketitle

\section{Introduction}
Fully nonlinear partial differential equations (PDEs) that incorporate elementary symmetric polynomials arise frequently in complex differential geometry.  Among these, complex $k$-Hessian equations are an important class of PDEs with applications in complex geometry, mathematical physics, and complex dynamical systems. Since Yau’s pioneering work that resolved the renowned Calabi conjecture in \cite{yau}, extensive research has been dedicated to investigating the existence, uniqueness, and fine properties of solutions to these equations, as well as developing innovative solution techniques. A notable example is the Fu--Yau equation \cite{FY07,FY08} arising from the Strominger system in complex geometry. In higher dimensions, the Fu--Yau equation corresponds precisely to an equation of complex $2$-Hessian type and has been extensively studied by Chu--Huang--Zhu \cite{CHZ19a,CHZ19b} and by Phong--Picard--Zhang \cite{PPZ16,PPZ17,PPZ19,PPZ21}. Another representative example  is the quaternionic Calabi conjecture on hyperK\"ahler with torsion (HKT) manifolds, proposed by Alesker--Verbitsky \cite{AV10}. In their work, the quaternionic Monge-Amp\`ere equation is reformulated as a special case of the complex Hessian equation on complex manifolds.

For the complex $k$-Hessian equations on compact manifolds, the central technical challenge lies in establishing a priori gradient estimates and second-order estimates for solutions. Under the non-negativity assumption on the holomorphic bisectional curvature, Hou \cite{Hou}, Jbilou \cite{Jbilou}, and Kokarev \cite{Kokarev} successfully derived the global a priori gradient estimates for solutions to these equations. 
Hou--Ma--Wu \cite{HMW} obtained a second-order estimate with quadratic dependence on the gradient. Then, Dinew--Ko\l odziej \cite{DK17} combined
 this estimate with a blow-up argument and a Liouville-type theorem to establish the gradient estimate.
On compact almost Hermitian manifolds, Chu--Huang--Zhang \cite{CHZ23}
 established a real Hessian estimate with quadratic gradient dependence for
 a general class of fully nonlinear elliptic equations admitting a
 $\mathcal C$-subsolution, including the complex Hessian equation.
For equations that also contain linear gradient terms,
 Huang--Zhang \cite{HZ} obtained a full real $C^2$ estimate
 with quadratic gradient dependence under a $\mathcal C$-subsolution
 assumption.

Very recently, Chu--Liu \cite{CL} have established a new gradient estimate that integrates a logarithmic modulus of continuity, the Hou--Ma--Wu estimate, and local $W^{2,p}$ estimates. Moreover, they explicitly clarified the precise dependence on the density. Their proof is built upon the stability property and on the auxiliary Monge--Amp\`ere equations,  instead of relying on blow-up and Liouville-type theorems. However, their proof still relies heavily on the Hou--Ma--Wu estimate. 
Guo--Song \cite{GS} have obtained gradient estimates for a broader class of concave elliptic equations on compact Hermitian manifolds under a fixed  smooth admissible subsolution with a uniform strict sup-slope gap.  Their proof uses stability and comparison with nearby smooth  admissible functions, without the Hou--Ma--Wu estimate or a  Liouville-type theorem. Under additional structural assumptions they also proved a separate Hessian estimate using Dirichlet Green functions. Their gradient theorem requires weaker assumptions on the density than those in \cite{CL}.

In Yau's proof of the Calabi conjecture \cite{yau}, the second-order estimate for the complex Monge--Amp\`ere equation on a compact K\"ahler manifold is obtained from a $C^0$ bound and the prescribed data, independently of any a priori gradient bound. The gradient estimate then follows from standard  interpolation. This suggests a similar strategy for complex Hessian equations: first establish a suitable second-order estimate, and then recover the gradient bound by interpolation. It is therefore natural to ask:
\begin{equation}\label{question}
\begin{gathered}
\text{\itshape Can one obtain a second-order estimate for complex Hessian equations}\\
\text{\itshape that yields a gradient bound by interpolation?}
\end{gathered}
\end{equation}
An estimate with linear dependence on the gradient is sufficient for this purpose:
\[
\sup_M\Delta_\omega u\leqslant  C\bigl(1+\sup_M|\nabla u|_\omega\bigr).
\]
Together with a $C^0$ bound, this inequality allows the gradient term to be absorbed by interpolation. Zhang--Zhang \cite{ZZ26} formulated the question in this linear form and proved the estimate under additional assumptions. In particular, they treated a regime with $2\leqslant  k\leqslant  n-1$ in which the wedge density is sufficiently close to $1$ in $L^1(M)$.

In the present work, we prove such a linear second-order estimate by the maximum principle in the nondegenerate K\"ahler setting. Interpolation then recovers the known gradient and complex Hessian bounds. This provides an alternative proof through the second-order estimate, without using an a priori gradient estimate, the Hou--Ma--Wu estimate, or Green functions. We now state our results more precisely.

Let $(M,\omega)$ be a compact K\"ahler manifold of complex
 dimension $n\ge2$. For a smooth real-valued function $u$, we write
\[
  \omega_u=\omega+\ddbar u,\qquad A_u=\omega^{-1}\omega_u.
\]
Let $\lambda(A)$ denote the eigenvalues of the Hermitian endomorphism $A$. For every fixed integer $2\leqslant  k\leqslant  n$,  the $k$-Hessian operator is given by
\[
 f(\lambda)=\sigma_k(\lambda)^{1/k}
\]
in the convex cone
\[
 \Gamma_k=\{\lambda\in\mathbb R^n:
\sigma_j(\lambda)>0\text{ for }1\leqslant  j\leqslant  k\},
\]
where the $j$-th elementary symmetric polynomial $\sigma_{j}$ is defined by
\[
\sigma_{j}(\lambda_1,\lambda_{2}, \cdots, \lambda_{n})=\sum_{1\leqslant   i_{1}<i_{2} <\cdots\mathrel{<}i_{j}\leqslant   n}\lambda_{i_1}\lambda_{i_2}\cdots \lambda_{i_{j}}.
\]
A function $u$ is called admissible if $\lambda(A_u)\in\Gamma_k$ on $M$. Denote
\[
F(A)=f(\lambda(A)).
\]
In the present work, we consider the nondegenerate complex $k$-Hessian equation
\begin{equation}\label{eq:equation}
  F(A_u)=h>0,\qquad \sup_M u=0, \qquad  \lambda(A_u)\in\Gamma_k.
\end{equation}
Equivalently,
\begin{equation}\label{eq:wedge}
  \omega_u^k\wedge\omega^{n-k}=\psi\omega^n,
  \qquad \psi=\frac{h^k}{\binom{n}{k}}.
\end{equation}
The compatibility condition $\int_M \psi\omega^n=\int_M\omega^n$ is automatic
 for an existing solution.
 
In local holomorphic coordinates, we write
$\omega=\sqrt{-1}\,g_{i\bar j}\,dz^i\wedge d\bar z^j$.
The inverse metric components $g^{i\bar j}$ are characterized by
\[
  \sum_{j=1}^n g^{i\bar j}g_{k\bar j}=\delta^i_k.
\]

Throughout the paper we use the conventions
\[
  \Delta_\omega u=g^{i\bar j}u_{i\bar j},\qquad
  |\nabla u|_\omega^2=g^{i\bar j}u_i u_{\bar j}.
\]
For convenience, we write $\lambda(A_{u})=(\lambda_1(A_{u}),\ldots,\lambda_n(A_{u}))$, with the eigenvalues arranged in decreasing order; that is,
\begin{equation}
\lambda_{1}(A_{u})\geqslant  \lambda_{2}(A_{u})\geqslant  \cdots \geqslant  \lambda_{n}(A_{u}).
\end{equation}
Therefore, 
\begin{equation*}
\tr A_u\equiv \sum_{i}\lambda_{i}(A_{u})=n+\Delta_\omega u.
\end{equation*}

Our main result is the following estimate.

\begin{theorem}\label{thm:main}
Fix positive constants $H_{1}$, $H_{2}$, $U$. Let $u$ be a smooth admissible solution
 of \eqref{eq:equation} such that
\begin{equation}\label{eq:data}
  h\geqslant H_{1}>0,\qquad \norm{h}_{C^2(M)}\leqslant  H_{2},
  \qquad \osc_M u\leqslant  U.
\end{equation}
There exists a constant $C$, depending only on
\[
  \data=(M,\omega,n,k,H_{1},H_{2},U),
\]
such that
\begin{equation}\label{eq:main}
  \sup_M\lambda_1(A_u)
\leqslant  C\left(1+\sup_M|\nabla u|_\omega^2\right)^{1/2}.
\end{equation}
Consequently,
\begin{equation}\label{eq:laplacian}
  \sup_M\Delta_\omega u
\leqslant  C\bigl(1+\sup_M|\nabla u|_\omega\bigr).
\end{equation}
\end{theorem}

The linear dependence in \eqref{eq:laplacian} allows us to close the
 gradient estimate by interpolation.

\begin{corollary}\label{cor:uniform}
Under the assumptions of Theorem \ref{thm:main}, there is a constant
 $C=C(\data)$ such that
\begin{equation}\label{eq:uniform}
  \norm{\nabla u}_{L^\infty(M)}
+\norm{\partial\bar\partial u}_{L^\infty(M)}\leqslant  C.
\end{equation}
\end{corollary}

We now briefly outline the core technical ingredients of the present work:
\begin{itemize}
\item[(1)]  In the prior work \cite{ZZ26}, the authors  omitted the crucial third-order term 
\begin{equation}\label{good third order term}
F^{i\overline{i}}\sum_{p>1}\frac{|u_{p\overline 1i}|^{2}+|u_{1\overline{p}i}|^{2}}{\lambda_{1}-\lambda_{p}}
\end{equation}
arising from the computation of $\mathcal{L}(\lambda_{1})$. In the present work, we invoke an improved lower bound for $\mathcal{L}(\log \lambda_{1})$ established by Guo--Song \cite{GS}, which makes essential use of the term \eqref{good third order term} and concavity.
\item[(2)] In order to derive \eqref{eq:main}, we employ a comparison argument inspired by the work of Guo--Song in \cite{GS}. More precisely, we choose a smooth admissible function $v$ close to $u$ such that $\norm{u-v}_{L^\infty}\leqslant  \varepsilon$ for a small constant $\varepsilon>0$. Then, we consider a new auxiliary function of the form
\begin{equation}\label{Q}
Q= \log\lambda_1
  -\frac12\log\left(1-\frac{|\nabla u|_\omega^2}{2K}\right)
  -Bu+N(u-v-\varepsilon)^2,
\end{equation}
where $K=1+\sup_M|\nabla u|_\omega^2$, $N=16B^2$, and $B>1$ is a sufficiently large constant with $\varepsilon=\frac{1}{64B}$. With this construction of the auxiliary function $Q$, we are able to establish a comparison estimate for second derivatives that exhibits linear dependence on both the gradient bound and $\norm{\nabla v}_{L^\infty}$.
\end{itemize}
The remaining ingredient is the $C^0(M)$ precompactness of the solution class, which follows from \cite{DK14,CL} independently of gradient and second-order estimates. Following the argument in the proof of \cite[Theorem 1.1]{GS}, we use a finite covering of the solution class in $C^0(M)$ to obtain the uniform estimate \eqref{eq:main} with a global constant independent of the individual solutions in $\mathcal{S}_{\mathcal{D}}$. 

At the end of the section, we give some further discussion on the auxiliary function $Q$ defined in \eqref{Q}.  Define a function $\gamma (t)=Bt-N(t-v-\varepsilon)^{2}$ on  a suitable interval, so that $Q$ can be written as
\[
Q= \log\lambda_1
  -\frac12\log\left(1-\frac{|\nabla u|_\omega^2}{2K}\right)-\gamma (u).
\]
A direct computation gives
\[
\gamma'(u)=B-2N(u-v-\varepsilon), \qquad \gamma''(u)=-2N.
\]
With the choices of $N=16B^2$  and $\varepsilon=\frac{1}{64B}$, it is not hard to verify the following important properties:
\[
\gamma'(u)\geqslant  B>1, \quad -\left[ \gamma''(u)+(\gamma'(u))^{2} \right] \geqslant  28B^2.
\]
As noted by B\l{}ocki in \cite{Blo09}, such a function $\gamma$ cannot be constructed  without introducing the nearby smooth admissible function $v$.

The paper is organized as follows. In Section \ref{sec:prelim}, we first collect some elementary inequalities about the $k$-Hessian operator, then derive a concavity inequality motivated by recent work of Dong--Zhang \cite{DZ},  and lastly prove compactness of the set of normalized smooth admissible solutions to \eqref{eq:equation}. In Section \ref{lower bound}, we prove an improved lower bound for $\mathcal{L}(\log \lambda_{1})$, by adapting an argument similar to that in \cite[Lemma 6.1]{GS}. In
Section \ref{sec:comparison}, we prove a comparison estimate for second derivatives. Finally, we complete the proofs of Theorem \ref{thm:main} and Corollary \ref{cor:uniform} in  Section \ref{sec:uniform}.

\section{Preliminaries}\label{sec:prelim}

\subsection{Elementary inequalities} For $f=\sigma_{k}^{1/k}$, we write
\[
  f_i=\frac{\partial f}{\partial\lambda_i},\qquad
  f_{ij}=\frac{\partial^2f}{\partial\lambda_i\partial\lambda_j},
  \qquad \mathcal{F}=\sum_i f_i.
\]
By homogeneity,
\begin{equation}\label{eq:linearized-u}
\sum_{i}f_{i}\lambda_{i}=f.
\end{equation}
The Newton--Maclaurin inequalities imply
\begin{equation}\label{eq:algebra}
  \mathcal{F}\geqslant \tau_{n,k}>0,\qquad
  \sigma_1(\lambda)\leqslant  C_{n,k}h\mathcal{F}^{k-1},
\end{equation}
for some positive constants $\tau_{n,k}$, $C_{n,k}$ depending only on $n$ and $k$. This is because
\[
  \mathcal{F}=\frac{n-k+1}{k}\,
\sigma_k^{1/k-1}\sigma_{k-1},\qquad
  \sigma_{k-1}^{k-1}\geqslant c_{n,k}\sigma_1\sigma_k^{k-2}.
\]
In particular, one can take $\tau_{n,k}=\binom{n}{k}^{1/k}$ in \eqref{eq:algebra}. Indeed, the concavity of $f$ at
 $\lambda$, applied to $\mathbf 1=(1,\ldots,1)$, and homogeneity give
\[
f(\mathbf 1)\leqslant  f(\lambda)
+\sum_if_i(1-\lambda_i)=\mathcal{F}.
\]
The ordering of the eigenvalues gives
\begin{equation}\label{eq:ordering}
  0<f_1\leqslant  \cdots\leqslant  f_n,\qquad f_n\geqslant \mathcal{F}/n.
\end{equation}
Since $\Gamma_k\subset\Gamma_2$,
\begin{equation}\label{eq:cone}
  \sum_i\lambda_i^2<\sigma_1(\lambda)^2,
  \qquad 0<\lambda_1\leqslant  \sigma_1(\lambda)\leqslant  n\lambda_1.
\end{equation}
These inequalities will also control the curvature terms in the
 covariant derivative calculations.

\subsection{A refined concavity inequality}
We now give an algebraic fact that we will use in the next section. This is just a restatement of
 \cite[Lemma 3.2]{DZ} adapted to our choice $f=\sigma_k^{1/k}$.

\begin{lemma}\label{lem:concavity}
For each $\beta\in(0,1/2)$ there exist a small constant $0<\delta_\beta\leqslant  \beta $ and a large constant
 $1\leqslant   R_\beta< \infty$, depending only on $n,k,H_{2},\beta$, such that if $\lambda \in\Gamma_k$ satisfies
 \begin{equation}
0< f(\lambda)\leqslant  H_{2}
 \end{equation}
and
 \begin{equation}\label{lambda 1n}
\lambda_1\geqslant R_\beta, \quad \lambda_n\geqslant -\delta_\beta\lambda_1,
 \end{equation}
 then, for every $\eta\in\mathbb C^n$,
\begin{equation}\label{eq:concavity}
\begin{split}
  -\sum_{p,q}f_{pq}\eta_p\overline{\eta_q}
  +\sum_{i>1}\frac{f_i|\eta_i|^2}{(1+\delta_\beta)\lambda_1}
  \ge{}&(1-\beta)\frac{f_1|\eta_1|^2}{\lambda_1}-\frac{k+1}{f}\left|\sum_if_i\eta_i\right|^2.
\end{split}
\end{equation}
\end{lemma}

\begin{proof}
The inequality in \cite[Lemma 3.2]{DZ} gives
\begin{equation}\label{DongZhang}
 -\frac{\sum_{p,q}\sigma_{k-2;pq}\eta_p\overline{\eta_q}}{\sigma_{k}}
 +2\frac{|\sum_i\sigma_{k-1;i}\eta_i|^2}{\sigma_{k}^2}
 +\sum_{i>1}\frac{\sigma_{k-1;i}|\eta_i|^2}{(1+\delta_\beta)\lambda_1\sigma_{k}}
 \geqslant (1-\beta)\frac{\sigma_{k-1;1}|\eta_1|^2}{\lambda_1\sigma_{k}},
\end{equation}
where
\[
\sigma_{k-1;i}=\frac{\partial \sigma_{k}}{\partial \lambda_{i}}, \quad \sigma_{k-2;ij}=\frac{\partial^{2}\sigma_{k}}{\partial \lambda_{i}\partial\lambda_{j}}.
\]
Multiplication by $f/k$, together with
\[
 f_p=\frac{f}{k\sigma_{k}}\sigma_{k-1;p},\qquad
 f_{pq}=\frac{f}{k\sigma_{k}}\sigma_{k-2;pq}
-\frac{(k-1)f}{k^2\sigma_{k}^2}\sigma_{k-1;p}\sigma_{k-1;q},
\]
gives the desired inequality \eqref{eq:concavity}.

To make these thresholds $\{R_{\beta},\delta_{\beta}\}$ uniform, we first apply  \cite[Lemma 3.2]{DZ}  at the
 fixed level $f(\mu)=1$, which yields the thresholds $\{\delta_{\beta},r_\beta\}$  that depend only on $n$, $k$ and $\beta$.
 
Then for a general $\lambda$, we define
\[
\mu=\frac{\lambda}{f(\lambda)}.
\]
The ratio $\frac{\lambda_{n}}{\lambda_1}$ remains invariant, i.e., $\frac{\lambda_{n}}{\lambda_1}=\frac{\mu_{n}}{\mu_{1}}$. Moreover, for any fixed $\eta\in\mathbb{C}^{n}$, each term in
 \eqref{eq:concavity} is homogeneous of degree $-1$
 in $\lambda$. As a consequence, the unit-level inequality formulated in terms of $\mu$ implies the desired inequality in terms of $\lambda$, provided that 
 \[
\mu_{1} \geqslant r_\beta, \qquad    \mu_{n}\geqslant -\delta_\beta \mu_{1}.
 \]
 These conditions reduce to $\lambda_{1}\geqslant  r_\beta f(\lambda)$ and $\lambda_{n}\geqslant  -\delta_{\beta}\lambda_{1}$.
It therefore suffices to set $R_\beta=\max\{1,H_{2}r_\beta\}$. This completes the proof.
\end{proof}

\subsection{Uniform continuity of the solutions}
In the lemma below, we prove a compactness result that does not depend  on  any norm of the derivatives of $u$, which serves as the crucial ingredient for establishing the Hessian estimate \eqref{eq:main}.
\begin{lemma}\label{lem:compactness}
For a fixed $\data$, the family of normalized smooth admissible solutions satisfying
 \eqref{eq:equation} and \eqref{eq:data}  is precompact
 in $C^0(M)$.
\end{lemma}

\begin{proof}
From the assumptions \eqref{eq:data}, we know that both the wedge densities $\psi$ in \eqref{eq:wedge} and the polynomial
 densities $\sigma_k(\lambda(A_u))=h^k$ are uniformly bounded in $L^p(M)$ for all exponents  $p>n/k$.  
 \begin{itemize}
\item[(i)]  For $1<k<n$, we apply the equicontinuity result of Dinew--Ko\l odziej \cite[Corollary 3.2]{DK14} to this family. \smallskip
\item[(ii)] For $k=n$, equicontinuity follows from the modulus of continuity in \cite[Theorem 3.1]{CL}. The same theorem also applies for each $1\leqslant  k\leqslant  n$.
\end{itemize}
The normalization $\sup_M u=0$ and the bound $\osc_M u\leqslant  U$ give $-U\leqslant  u\le0$.
The Arzel\`a--Ascoli theorem therefore gives precompactness in $C^0(M)$.
This compactness property follows from zeroth order stability alone,
 independently of any gradient or second-order estimates.
\end{proof}

\begin{remark}
It should be noted that the derivation of the modulus of continuity in \cite[Section 3]{CL} does not make use of the Hou--Ma--Wu estimate, which only appears in their subsequent gradient estimates. For our argument, only compactness is required; no quantitative modulus of continuity appears in  the  maximum principle argument that follows.
\end{remark}

\section{Derivatives of the largest eigenvalue}\label{lower bound}
Let $\nabla$ be the Chern connection of $\omega$, which agrees with the
Levi--Civita connection in the K\"ahler setting. Derivative indices are appended
in the order in which the covariant derivatives are taken:
\[
  T_{;i\bar j}=\nabla_{\bar j}\nabla_iT,\qquad
  u_{p\bar q i\bar j}=\nabla_{\bar j}\nabla_i(u_{p\bar q}).
\]
We omit semicolons in derivative subscripts and use the curvature convention
\[
  \Gamma^r_{ip}=g^{r\bar s}\partial_i g_{p\bar s},\qquad
  R_{i\bar j p\bar q}=-g_{r\bar q}\partial_{\bar j}\Gamma^r_{ip}.
\]
In normal holomorphic coordinates, this gives
$R_{i\bar j p\bar q}=-\partial_i\partial_{\bar j}g_{p\bar q}$ at the center.
Metric compatibility gives $\nabla g=0$, so all covariant derivatives of $g$
vanish. Lowering the endomorphism index of $A_u$ with $g$, we write
\[
  A^u_{p\bar q}=(\omega_u)_{p\bar q}=g_{p\bar q}+u_{p\bar q}.
\]
In an $\omega$-unitary frame, this Hermitian matrix represents $A_u$.
For the matrix derivatives of $F$, we write
\[
  F^{i\bar j}=\frac{\partial F}{\partial A^{u}_{i\bar j}},
  \qquad F^{i\bar j,k \bar l}=\frac{\partial^2 F}{\partial A^{u}_{i\bar j} \partial A^{u}_{k\bar l}},
  \qquad |\partial w|_F^2=F^{i\bar j}w_iw_{\bar j}.
\]
The pointwise component calculations below use an $\omega$-unitary frame
diagonalizing $A_u$. In this frame,
\[
  g_{i\bar j}=\delta_{ij},\qquad
  A^u_{i\bar j}=\lambda_i\delta_{ij},\qquad F^{i\bar j}=f_i\delta_{ij},
\]
with no summation over the free indices $i,j$ in these identities.
The linearized operator associated with \eqref{eq:equation} is given by
\[
 \Lin(w)=F^{i\bar j}\nabla_i\nabla_{\bar j}w.
\]
Notice that the assumption \eqref{lambda 1n} will be checked
 at the auxiliary maximum in one case of the proof. As a result, our use of
 \cite[Lemma 3.2]{DZ} imposes no additional global eigenvalue condition on the solution.

In this section, we will establish an improved lower bound for $\mathcal{L}$ applied to the logarithm of the largest eigenvalue.
\begin{lemma}\label{lem:eigenvalue}
Let $u$ be a smooth admissible solution of \eqref{eq:equation}. Then 
\begin{itemize}
\item[(1)] There exists a constant $C_{0}=C_{0}(\data)$ such that 
\begin{equation}\label{eq:rough}
  \Lin(\log\lambda_1)
     \geqslant -|\partial\log\lambda_1|_F^2-C_{0}\mathcal{F}.
\end{equation}
\item[(2)] Moreover, for every $\beta\in(0,1/2)$ such that $\lambda(A_{u})$ fulfills \eqref{lambda 1n}, we have
\begin{equation}\label{eq:refined}
  \Lin(\log\lambda_1)
     \geqslant -\beta|\partial\log\lambda_1|_F^2-C_{1}\mathcal{F}.
\end{equation}
Here the constant $C_{1}=C_{1}(\mathcal{D})$ is independent of $\beta$.
\end{itemize}
At multiple eigenvalues, both inequalities hold for smooth functions touching $\lambda_1$ from above, with $\lambda_1$ replaced by the test
 function on the left and in its gradient.
\end{lemma}

\begin{proof}
The proof originates from \cite[Lemma 6.1]{GS}, and we include the calculation here for the sake of completeness.

{\em Step 1.} First suppose that $\lambda_1$ is simple. In a normal unitary frame
 diagonalizing $A_u$, for each $p=1,\cdots,n$ we write
\begin{equation}
  \eta_p=\nabla_1(\omega_u)_{p\bar p},\qquad
  z_p=\nabla_p(\omega_u)_{1\bar1}.
\end{equation}
Closedness of $\omega_u$ gives the Codazzi identities, in particular
\begin{equation}
z_1=\eta_1, \qquad     z_p=\nabla_1(\omega_u)_{p\bar1}.
\end{equation}
Differentiating equation \eqref{eq:equation} once yields
\begin{equation}\label{eq:first-equation}
  \sum_pf_p\eta_p=h_1.
\end{equation}
At a diagonal matrix with simple largest eigenvalue $\lambda_1$, the standard derivative formulas are as follows (see \cite{Gerhardt,Spruck}):
\begin{equation}\label{derivatives of lambda}
\begin{split}
\lambda_{1}^{pq}& =\delta_{1p}\delta_{1q}\,,\\
\lambda_{1}^{pq,rs}&=
\begin{cases}
\dfrac{1}{\lambda_1-\lambda_p}, & q=r=1,\ p=s>1,\\[1mm]
\dfrac{1}{\lambda_1-\lambda_r}, & p=s=1,\ r=q>1,\\[1mm]
0, & \text{otherwise},
\end{cases}
\end{split}
\end{equation}
where $\lambda_{1}^{pq}$ represents the derivative with respect to the entry located at the $p$-th row and $q$-th column. 
It follows from \eqref{derivatives of lambda} along with the K\"ahlerity of $\omega$ that

\begin{equation}
\begin{split}
\lambda_{1,i}&=\lambda_{1}^{pq}(A^{u}_{p\bar q})_{i}=u_{1\overline{1}i}\,,\\
\lambda_{1,i\overline{i}}&=\lambda_{1}^{pq,rs}(A^{u}_{p\bar q})_{i}(A^{u}_{r\bar s})_{\overline i}+\lambda_{1}^{pq}(A^{u}_{p\bar q})_{i\overline{i}}\\
&=u_{1\overline{1}i\overline i}+\sum_{p>1}\frac{|u_{p\overline 1i}|^{2}+|u_{1\overline{p}i}|^{2}}{\lambda_{1}-\lambda_{p}}\,.
\end{split}
\end{equation}

By commuting the covariant derivatives, we have
\begin{equation}\label{nonnegative 1}
\begin{split}
u_{1\overline{1}i\overline i}&=u_{i\overline i1\overline{1}}+(u_{1\overline{1}}-u_{i\overline i})R_{1\overline 1i\overline{i}}\\
&=u_{i\overline i1\overline{1}}+(\lambda_{1}-\lambda_{i})R_{1\overline 1i\overline{i}}\\
&\geqslant   u_{i\overline i1\overline{1}}+(\lambda_{1}-\lambda_{i})\inf_{l\neq 1}R_{1\overline 1l\overline{l}}\,.
\end{split}
\end{equation}
It follows that

\begin{equation}\label{AA}
\begin{split}
\mathcal{L}(\lambda_{1})\geqslant   & \, F^{i\overline i}u_{i\overline i1\overline{1}}\mathbin{+}\inf_{l\neq 1}R_{1\overline 1l\overline{l}}\sum_{i=1}^{n}F^{i\overline i}(\lambda_{1}-\lambda_{i})+F^{i\overline{i}}\sum_{p>1}\frac{|u_{p\overline 1i}|^{2}+|u_{1\overline{p}i}|^{2}}{\lambda_{1}-\lambda_{p}}\\
\geqslant  &\,  F^{i\overline i}u_{i\overline i1\overline{1}}-C_{R}\sum_{i=1}^{n}F^{i\overline i}(\lambda_{1}-\lambda_{i})+\sum_{i>1}\frac{F^{i\overline{i}}|u_{1\overline{i}i}|^{2}}{\lambda_{1}-\lambda_{i}}+\sum_{p>1}\frac{F^{1\overline{1}}|u_{p\overline 11}|^{2}}{\lambda_{1}-\lambda_{p}}\\
\geqslant  &\,  F^{i\overline i}u_{i\overline i1\overline{1}}-C_{R}\sum_{i=1}^{n}F^{i\overline i}(\lambda_{1}-\lambda_{i})+\sum_{i>1}\frac{f_{i}|\eta_{i}|^{2}}{\lambda_{1}-\lambda_{i}}+\sum_{i>1}\frac{f_{1}|z_{i}|^{2}}{\lambda_{1}-\lambda_{i}}\,,
\end{split}
\end{equation}
where $C_{R}$ is a constant depending on the holomorphic bisectional curvature.
Differentiating equation \eqref{eq:equation} twice, we see that

\begin{equation*}
F^{i\overline i}u_{i\overline{i}1\overline{1}}+F^{i\overline{j},k\overline{l}}u_{i\overline{j}1} u_{k\overline{l}\overline{1}}=h_{1\overline{1}}\,.
\end{equation*}

Inserting it into \eqref{AA}, we infer that

\begin{equation*} 
\begin{split}
\mathcal{L}(\lambda_{1})\geqslant   &-F^{i\overline{j},k\overline{l}}u_{i\overline{j}1} u_{k\overline{l}\overline{1}}-C_{R}\sum_{i=1}^{n}F^{i\overline i}(\lambda_{1}-\lambda_{i})+h_{1\overline{1}}\\
&+\sum_{i>1}\frac{f_{i}|\eta_{i}|^{2}}{\lambda_{1}-\lambda_{i}}+\sum_{i>1}\frac{f_{1}|z_{i}|^{2}}{\lambda_{1}-\lambda_{i}}\\
\geqslant   & -F^{i\overline{j},k\overline{l}}u_{i\overline{j}1} u_{k\overline{l}\overline{1}}-C_{R}(1+\lambda_{1})\mathcal{F}+h_{1\overline{1}}+\sum_{i>1}\frac{f_{i}|\eta_{i}|^{2}}{\lambda_{1}-\lambda_{i}}+\sum_{i>1}\frac{f_{1}|z_{i}|^{2}}{\lambda_{1}-\lambda_{i}}\,.
\end{split}
\end{equation*}

Here, we have used the homogeneity of $F$ in the last inequality. Notice that
\begin{equation*}
F^{i\overline{j},k\overline{l}}=
\left\{\begin{array}{ll}
F^{i\overline{i},k\overline{k}}, \quad &\textrm{ if }  i=j, k=l,  \\[1mm]
F^{i\overline{k},k\overline{i}}, \quad &\textrm{if } i=l,\ k=j,\ i\ne k,\\[1mm]
0, \quad & \textrm{ otherwise.}
\end{array}\right.
\end{equation*}
Therefore, by concavity of $F$,
\begin{equation}
\begin{split}
-F^{i\overline{j},k\overline{l}}u_{i\overline{j}1} u_{k\overline{l}\overline{1}}\geqslant  & -F^{i\overline{i},k\overline{k}}u_{i\overline{i}1} u_{k\overline{k}\overline{1}}\\
&\quad-\sum_{i>1}F^{i\overline{1},1\overline{i}}\left(|u_{i\overline{1}1}|^{2}+|u_{1\overline{i}1}|^{2}\right)\\
\geqslant  & -\sum_{p,q}f_{pq}\eta_p\overline{\eta_q}+\sum_{i>1}\frac{(f_{i}-f_{1})|z_{i}|^{2}}{\lambda_{1}-\lambda_{i}}.
\end{split}
\end{equation}
Consequently,  we have
\begin{equation}\label{eq:spectral}
\begin{split}
\mathcal{L}(\lambda_{1})\ge{}&h_{1\bar1}-C_{R}(1+\lambda_{1})\mathcal{F}
-\sum_{p,q}f_{pq}\eta_p\overline{\eta_q}
 +\sum_{i>1}\frac{f_i|\eta_i|^2}{\lambda_1-\lambda_i}
+\sum_{i>1}\frac{f_i|z_i|^2}{\lambda_1-\lambda_i}.
\end{split}
\end{equation}
Discarding the nonnegative quadratic terms by concavity, dividing
 by $\lambda_1$, and using $\mathcal{F}\geqslant \tau_{n,k}$ proves \eqref{eq:rough}.

In the region of \eqref{eq:refined}, it suffices to prove that there exists a constant $C$, independent of $\beta$, such that
\begin{equation}\label{eq:refined two}
\begin{split}
\mathcal{L}(\lambda_{1}) & \geqslant  (1-\beta) \frac{|\partial \lambda_{1}|^{2}_{F}}{\lambda_{1}} -C\lambda_{1}\mathcal{F}\\
&= (1-\beta)\frac{f_{1}|\eta_{1}|^{2}}{\lambda_{1}}+(1-\beta) \sum_{i>1}\frac{f_{i}|z_{i}|^{2}}{\lambda_{1}}-C\lambda_{1}\mathcal{F}
\end{split}
\end{equation}
 provided that 
\begin{equation}\label{lambda 1i}
\lambda_1-\lambda_i\le(1+\delta_\beta)\lambda_1  \qquad (i> 1).
\end{equation}
On the one hand,  using \eqref{eq:concavity}, \eqref{eq:first-equation} and \eqref{lambda 1i}, we have
\begin{equation}\label{third 1}
\begin{split}
-\sum_{p,q}f_{pq}\eta_p\overline{\eta_q}
 +\sum_{i>1}\frac{f_i|\eta_i|^2}{\lambda_1-\lambda_i} & \geqslant  -\sum_{p,q}f_{pq}\eta_p\overline{\eta_q}
 +\sum_{i>1}\frac{f_i|\eta_i|^2}{(1+\delta_\beta)\lambda_1} \\
 & \geqslant  (1-\beta)\frac{f_1|\eta_1|^2}{\lambda_1}-(k+1)\frac{|h_1|^2}{h}.
\end{split}
\end{equation}
On the other hand, since $\delta_\beta\leqslant  \beta$, we have $\frac{1}{1+\delta_{\beta}}\geqslant  1-\beta$. Hence, it follows from \eqref{lambda 1i} that
\begin{equation}\label{third 2}
\sum_{i>1}\frac{f_i|z_i|^2}{\lambda_1-\lambda_i}\geqslant   (1-\beta) \sum_{i>1}\frac{f_{i}|z_{i}|^{2}}{\lambda_{1}}.
\end{equation}
Putting \eqref{third 1} and \eqref{third 2} into \eqref{eq:spectral}, we infer that 
\begin{equation}
\begin{split}
\mathcal{L}(\lambda_{1})\ge{}&h_{1\bar1}-(k+1)\frac{|h_1|^2}{h}-C_{R}(1+\lambda_1)\mathcal{F}+
(1-\beta)\frac{f_{1}|\eta_{1}|^{2}}{\lambda_{1}}+(1-\beta) \sum_{i>1}\frac{f_{i}|z_{i}|^{2}}{\lambda_{1}}.
\end{split}
\end{equation}
The term involving the prescribed density
 is bounded below independently of $\beta$:
\[
  h_{1\bar1}-(k+1)\frac{|h_1|^2}{h}\geqslant -C_{2}
\]
for some positive constant $C_{2}$ depending only on $H_{1}$, $H_{2}$ and $k$. Since $\mathcal{F}\geqslant \tau_{n,k}$ and $\lambda_{1}\geqslant  1$,
 this proves \eqref{eq:refined two} with $C$ independent of $\beta$.\medskip

{\em Step 2.} Now suppose that $\lambda_{1}$ has multiplicity $m>1$, and that
 a smooth function $\ell$ touches $\lambda_1$ from above.
First order contact implies
\[
  \nabla_i(\omega_u)_{p\bar q}=\ell_i\delta_{pq}
\qquad(p,q\leqslant  m).
\]
For $i\leqslant  m$, choosing $p\leqslant  m$ with $p\ne i$ and applying the Codazzi identity
 shows $\ell_i=\nabla_p(\omega_u)_{i\bar p}=0$.
It follows that $\eta_p=0$ for $p\leqslant  m$ and $z_i=0$ for $i\leqslant  m$.
The upper-contact second derivative formula retains the variation
 terms with positive gaps, $i>m$. All the omitted numerators vanish,
 so the same calculation applies to $\ell$.
\end{proof}

\section{The comparison estimate}\label{sec:comparison}

We first prove a comparison estimate in which the permitted $C^0$
 distance is independent of the comparison function.

\begin{proposition}\label{prop:comparison}
There exists a small constant $\varepsilon=\varepsilon(\data)>0$ such that the following holds.
Let $u$ satisfy the assumptions of Theorem \ref{thm:main}, and let $v$
 be any smooth admissible function with
\[
  \norm{u-v}_{L^\infty(M)}\leqslant  \varepsilon.
\]
Writing $V=\norm{\nabla v}_{L^\infty(M)}$, one has
\begin{equation}\label{eq:comparison-estimate}
  \sup_M\lambda_1(A_u)
\leqslant  C(\data,V)\left(1+\sup_M|\nabla u|_\omega^2\right)^{1/2}.
\end{equation}
\begin{remark}
We remark that during the proof of \eqref{eq:comparison-estimate}, no second derivative bound for $v$ is required. Instead, we use only the fact that $\lambda(A_v)\in\Gamma_k$.
\end{remark}
\end{proposition}

\begin{proof}
Set
\[
  P=|\nabla u|_\omega^2,\qquad K=1+\sup_MP,\qquad
  \phi(s)=-\frac12\log\left(1-\frac{s}{2K}\right).
\]
For $0\leqslant  s\leqslant  K$, we have
\begin{equation}\label{eq:phi}
  \frac1{4K}\leqslant  \phi'\leqslant  \frac1{2K},\qquad
  \phi''=2(\phi')^2,\qquad 0\leqslant  \phi(P)\leqslant  \frac12\log2.
\end{equation}
Choose a sufficiently large constant $B\geqslant 1$, depending only on
 $\data$, and fix
\begin{equation}\label{eq:parameters}
  N=16B^2,\qquad \varepsilon=\frac1{64B}.
\end{equation}
The required lower bound for $B$ will be specified below and will
 not involve $v$. Given $v$ as in the proposition, define
\[
  w=u-v-\varepsilon\in  [-2\varepsilon, 0].
\]
Consider the auxiliary function
\begin{equation}\label{eq:test}
  Q=\log\lambda_1+\phi(P)-Bu+Nw^2.
\end{equation}

We first calculate the terms in $\Lin (Q)$.
Concavity of $F$ and admissibility of $v$ imply
\[
  \Lin (w)\leqslant  F(A_u)-F(A_v)\leqslant  h.
\]
Since $w\le0$ and $h\leqslant   H_{2}$, it follows that
\begin{equation}\label{eq:q-square}
  \Lin(Nw^2)=2N w\Lin (w)+2N|\partial w|_F^2
     \geqslant -4N\varepsilon H_{2}+2N|\partial w|_F^2.
\end{equation}
This is the point at which admissibility of $v$ replaces a bound for its second derivatives.

A direct calculation shows that
\begin{equation*} 
  \Lin (P)=   F^{i\overline{i}}\Big(|u_{i\overline i}|^{2}+\sum_{p=1}^{n}|u_{ip}|^{2}\Big)+\sum_{p=1}^{n}F^{i\overline i}(u_{pi\overline i}u_{\overline p}+u_{\overline pi\overline i}u_{p})\,.
\end{equation*}
By commuting the covariant derivatives,
\[
u_{pi\overline i}=u_{i\overline ip}+\sum_{q=1}^{n}u_{q}R_{i\overline ip\overline{q}}\,, \qquad u_{\overline pi\overline i}=u_{i\overline i\overline p}\,.
\]
Differentiating \eqref{eq:equation} once, we obtain
\begin{equation*}
F^{i\overline i}u_{i\overline{i}p}=h_{p}\,, \qquad F^{i\overline i}u_{i\overline{i}\overline p}=h_{\overline p}\,.
\end{equation*}
Thus, we infer that
\begin{equation}\label{nonnegative 2}
\begin{split}
\sum_{p=1}^{n}F^{i\overline i}(u_{pi\overline i}u_{\overline p}+u_{\overline pi\overline i}u_{p})
&= 2\sum_{p=1}^{n}\mathrm{Re}\Big(h_{p}u_{\overline{p}} \Big)+F^{i\overline{i}}\sum_{p,q=1}^{n}u_{\overline{p}}u_{q}R_{i\overline ip\overline{q}}\\
&\geqslant   -C_{3}\sqrt{P}-C_{R}P\mathcal{F}
\end{split}
\end{equation}
for some constant $C_{3}$ depending only on $\|\nabla h\|_{L^{\infty}}$. It follows from \eqref{nonnegative 2} that
\begin{equation}\label{eq:bochner}
  \Lin (P)\geqslant \sum_if_i(\lambda_i-1)^2-C_{3}\sqrt{P}-C_{R}P\mathcal{F}.
\end{equation}
Using the elementary inequality $(\lambda_i-1)^2\geqslant \frac{1}{2}\lambda_i^2-1$ and \eqref{eq:phi}, we get
\begin{equation}\label{eq:gradient-positive}
  \phi'\Lin (P)\geqslant \frac1{8K}\sum_if_i\lambda_i^2-C_{R}\mathcal{F}-C_{3}.
\end{equation}
In addition, \eqref{eq:linearized-u} implies
\begin{equation}\label{eq:linear-weight}
  \Lin(-Bu)=B\mathcal{F}-Bh.
\end{equation}

Let $x_0$ be a maximum point of $Q$. We consider two cases.
\begin{itemize}
\item[(i)] Assume $\lambda_{1}$ is simple at $x_{0}$; that is, $\lambda_1$ is smooth there. At $x_0$,
\begin{equation}\label{eq:first-order}
  \partial\log\lambda_1
      =-\phi'\partial P+(B-2N w)\partial w+B\partial v.
\end{equation}
Since $B-2N w\leqslant  B+4N\varepsilon=2B$, two applications of the elementary square inequality yield
\begin{equation}\label{eq:gradient-loss}
  |\partial\log\lambda_1|_F^2
   \leqslant  2(\phi')^2|\partial P|_F^2
+16B^2|\partial w|_F^2+4B^2V^2\mathcal{F}.
\end{equation}
\item[(ii)] Otherwise,  the function
\[
  \ell(x)=\exp\bigl(Q(x_0)-\phi(P(x))+Bu(x)-N w(x)^2\bigr)
\]
touches $\lambda_1$ from above at $x_0$.
Lemma \ref{lem:eigenvalue} applied to $\ell$ justifies the same  calculation. In particular, the corresponding smooth test for $Q$  is constant, which is sufficient for the maximum principle inequality.
\end{itemize}

Let $\tilde{\beta}\in(0,1]$ denote the loss coefficient in Lemma \ref{lem:eigenvalue}; that is,
\begin{equation}
\tilde{\beta}=
\left\{\begin{array}{ll}
\beta, \quad &\textrm{ if \eqref{lambda 1n} holds},  \\[1mm]
1,  \quad &\textrm{otherwise}.
\end{array}\right.
\end{equation}
The first two terms on the right of \eqref{eq:gradient-loss}, after
 multiplication by $\tilde{\beta}$, are absorbed because
\[
  \phi''-2\tilde{\beta}(\phi')^2\ge0,\qquad
  2N-16\tilde{\beta} B^2\geqslant 16B^2.
\]
Combining \eqref{eq:q-square}, \eqref{eq:gradient-positive},
 \eqref{eq:linear-weight} and Lemma \ref{lem:eigenvalue}, we obtain
 at $x_0$, whenever the relevant eigenvalue threshold is satisfied,
\begin{equation}\label{eq:key}
  0\geqslant \frac1{8K}\sum_if_i\lambda_i^2
+\bigl(B-C_{0}-C_{1}-C_{R}-4\tilde{\beta} B^2V^2\bigr)\mathcal{F}-C_{4}.
\end{equation}
Here $C_{0}$ and $C_{1}$ are the coefficients of $\mathcal{F}$ in Lemma \ref{lem:eigenvalue}. We may take
\[
  C_{4}\equiv C_{3}+(B+4N\varepsilon)H_{2}=C_{3}+2BH_{2},
\]
where $C_{3}$ is the constant term in \eqref{eq:gradient-positive}.
The constant $C_{3}$ depends only on $\data$ and is independent
 of $B,K,V,\tilde{\beta}$.
We now specify the earlier choice by taking
\begin{equation}\label{eq:a}
  B\equiv C_{0}+C_{1}+C_{R}+2.
\end{equation}
Thus $B,N,\varepsilon$ are determined solely by $\data$.

After $v$ has been fixed, choose
\begin{equation}\label{eq:beta}
  0<\beta\leqslant  \min\left\{\frac12,\,\frac1{8B^2(1+V^2)}\right\},
\end{equation}
and let $\delta_\beta,R_\beta$ be as in
 Lemma \ref{lem:eigenvalue}.
If $\lambda_1(x_0)<R_\beta$, the required bound at $x_0$ is immediate.
Otherwise, there are two cases.

\smallskip
\emph{Case 1: If $\lambda_n(x_0)\geqslant -\delta_\beta\lambda_1(x_0)$.}
In this case we treat \eqref{eq:key} with $\tilde{\beta}=\beta$.
By \eqref{eq:a} and \eqref{eq:beta}, we have
\[B-C_0-C_{1}-C_{R}-4\beta B^2V^2\geqslant 1\]
and hence $\mathcal{F}(x_0)\leqslant  C_{4}$.
It follows from \eqref{eq:algebra} that
\begin{equation}\label{eq:case-positive}
  \lambda_1(x_0)\leqslant  \sigma_1(\lambda(x_0))
\leqslant  C_{n,k}h(x_0)\mathcal{F}(x_0)^{k-1}\leqslant  C_{5}
\end{equation}
for a uniform  constant $C_5$ depending only on $\data$.

\smallskip
\emph{Case 2: If $\lambda_n(x_0)<-\delta_\beta\lambda_1(x_0)$.}
In this case we treat \eqref{eq:key} with $\tilde{\beta}=1$.
By \eqref{eq:ordering}, we have
\[
  \sum_if_i\lambda_i^2\geqslant f_n\lambda_n^2
\geqslant \frac{\delta_\beta^2}{n}\lambda_1^2\mathcal{F}.
\]
Consequently,
\[
  \frac{\delta_\beta^2}{8nK}\lambda_1(x_0)^2\mathcal{F}
      \leqslant  4B^2V^2 \mathcal{F}+C_{4}.
\]
Since $\mathcal{F}\geqslant \tau_{n,k}$, we obtain
\begin{equation}\label{eq:case-negative}
  \lambda_1(x_0)\leqslant  C(\data,V)\sqrt K.
\end{equation}
Here the dependence on $V$ includes that of
 $\beta,\delta_\beta,R_\beta$.

In every case, we have $\lambda_1(x_0)\leqslant  C(\data,V)\sqrt K$.
The oscillation of the remaining part of $Q$ is bounded by
\[
  \osc_M\bigl(\phi(P)-Bu+Nw^2\bigr)
      \leqslant  \tfrac12\log2+BU+4N\varepsilon^2.
\]
Comparing $Q(x)$ with $Q(x_0)$ proves
 \eqref{eq:comparison-estimate}.
\end{proof}

\section{Uniformity and the gradient bound}\label{sec:uniform}
\subsection{Proof of Theorem \ref{thm:main}}

Let $\mathcal S_{\data}$ denote the family of normalized smooth admissible solutions of \eqref{eq:equation} satisfying
 \eqref{eq:data}. The right-hand side may vary within the prescribed
 bounds. If this class is empty, the assertion is vacuous. Suppose that the desired estimate
\eqref{eq:main} does not hold. Then there exists a sequence $\{u_{j}\}_{j\in \mathbb{N}}\subset \mathcal S_{\data}$ such that
\[
\frac{\sup_M\lambda_1(A_{u_{j}})}{\left(1+\sup_M|\nabla u_{j}|_\omega^2\right)^{1/2}}\longrightarrow \infty.
\]

Fix $\varepsilon=\varepsilon(\data)$ from Proposition \ref{prop:comparison}.
By Lemma \ref{lem:compactness}, $\mathcal S_{\data}$ is totally bounded in $C^0(M)$.
It therefore admits a finite $\varepsilon$-net with centers in $\mathcal S_{\data}$, so we may choose finitely many
smooth admissible functions $v_1,\ldots,v_N\in\mathcal S_{\data}$ such that
\begin{equation}\label{eq:finite-cover}
  \mathcal S_{\data}\subset
  \bigcup_{\ell=1}^N
  \bigl\{w\in C^0(M):
 \norm{w-v_\ell}_{L^\infty(M)}<\varepsilon\bigr\}.
\end{equation}
The centers are chosen in the smooth solution class itself.
Fix this finite family once for the class determined by $\data$. Therefore, there exist an index $\ell_{0}\in \{1,2,\cdots,N\}$ and  a subsequence $\{u_{j_{k}}\}_{k\in\mathbb{N}}\subset \mathcal{S}_{\data}$ such that 
\[
\|u_{j_{k}}-v_{\ell_{0}}\|_{L^{\infty}}\leqslant  \varepsilon \qquad (\text{for all sufficiently large } k).
\]
Then Proposition \ref{prop:comparison} gives 
\begin{equation}\label{eq:uniform-constant}
\frac{\sup_M\lambda_1(A_{u_{j_{k}}})}{\left(1+\sup_M|\nabla u_{j_{k}}|_\omega^2\right)^{1/2}} \leqslant  C(\mathcal{D},V_{\ell_{0}}),
\end{equation}
where we denote $V_{\ell}=\|\nabla v_{\ell}\|_{L^{\infty}(M)}$. The right-hand side is finite and independent of $k$, which yields a contradiction. This completes the proof of Theorem \ref{thm:main}.

\subsection{Proof of Corollary \ref{cor:uniform}}
\begin{proof}
Without loss of generality, we may assume $G=\norm{\nabla u}_{L^\infty(M)}>1$. By compactness of $M$, there exists a point $p\in M$ such that \[G=|\nabla u|(p).\]
Fix an exponent $q>2n$. Local elliptic $W^{2,q}$-estimates and Sobolev embedding give constants $r_0\in(0,1)$ and $\overline{C}$, depending only on $(M,\omega)$ and the fixed exponent, such that for $0<r<r_0$,
\begin{equation*} 
\begin{split}
 r|\nabla u|(p)&\leqslant    \overline{C}\operatorname{osc}_{B_{r}(p)}u+\overline{C}r^{2}\norm{\Delta_\omega u}_{L^{\infty}({B_{r}(p)})}.
\end{split}
\end{equation*}
Here $B_{r}(p)\Subset M$ denotes the coordinate ball centered at $p$ with radius $r$. 
Notice that Theorem \ref{thm:main} yields
\[
\norm{\Delta_\omega u}_{L^\infty(M)}\leqslant  C(\mathcal{D})(1+G)\leqslant   2C(\mathcal{D}) G.
\]
Therefore,
\begin{equation}\label{rG}
 rG \leqslant   \overline{C} U+2\overline{C}C(\data)r^{2}G.
\end{equation}
If $G\leqslant  r_0^{-2}$, the gradient is already bounded. Otherwise, choose $r=G^{-1/2}<r_0$. Then \eqref{rG} gives
\[
G\leqslant   \left(\overline{C} U+2\overline{C}C(\data)\right)^{2}.
\]
Thus $G\leqslant  C(\data)$ in all cases. Theorem \ref{thm:main} now bounds $\lambda_1(A_u)$, and \eqref{eq:cone} bounds all the eigenvalues of $A_u$. Since $A_u=\mathrm{Id}+\omega^{-1}\ddbar u$, this also gives the asserted complex Hessian bound and completes the proof.
\end{proof}

Our argument does not require an explicit modulus of continuity in the range $1<k<n$, in contrast to the proof of Chu--Liu \cite[Theorem 1.2]{CL}. In this range, the equicontinuity result of Dinew--Ko\l{}odziej \cite{DK14} suffices for the precompactness used in Lemma \ref{lem:compactness}.

\subsection{Further questions}
We conclude with two questions concerning possible extensions of this proof.

\subsubsection{The degenerate setting}
Theorem \ref{thm:main} assumes a positive lower bound for $h$. Under what regularity assumptions on a nonnegative wedge density $\psi$ can we use the present argument to obtain \eqref{eq:main} uniformly along positive approximations? The distinction between the regularity of $\psi$ and that of its root $h$ is important: smoothness of a nonnegative $\psi$ need not give uniform $C^2$ control of $h$. On the other hand, uniform $C^2$ bounds for nonnegative root densities control $|\nabla h|^2/h$ on $\{h>0\}$ independently of a positive lower bound. A passage to degenerate solutions would also require tracking the constants and allowing admissibility in $\overline{\Gamma_k}$. For related degenerate estimates, see Chu--McCleerey \cite{CM21}.

\subsubsection{Non-K\"ahler setting}
Can we adapt the present maximum principle argument to compact Hermitian manifolds? The K\"ahler condition enters the Codazzi identities and the largest-eigenvalue calculation. On a Hermitian background, torsion introduces additional third-order terms, although the gradient Bochner terms can still be controlled; see \cite[Lemma 3.4]{GS}. The main issue for the present method is to control the torsion with error bounds independent of the small parameter $\beta$. This uniformity allows us to fix the comparison distance before choosing the nearby smooth admissible function. Further cancellations or a modified auxiliary function may be needed to preserve this order of choices; see \cite[Section 4]{DZ} for the torsion terms arising in Hermitian estimates.\medskip

\noindent\textbf{Acknowledgements}
The first author was supported by the NSFC (No. 12401075) and
Fundamental Research Funds for the Central Universities (Grant No. 20720260220).  The second author was supported by the NSFC (No. 12401065). The third author was supported by the National Key R and D Program of China 2020YFA0713100 and by NSFC (No. 12141104, 11801535  and 11721101). \medskip

\noindent\textbf{Use of AI.} This work is a continuation of arXiv:2610.07852. Compared with the previous version, the key improvement in the present work relies on the techniques of Guo and Song \cite{GS} for handling the concavity of the Hessian operators.  During the preparation of this work, we use AI tools, particularly
ChatGPT 6.1 Astra, to explore proof strategies and assist with calculations.
They verified and revised all AI-generated suggestions incorporated into
the paper and take full responsibility for its contents.

\end{document}